\documentclass[12pt]{article}

\usepackage{amsfonts,amsmath,amsthm,latexsym,amssymb,fullpage,amscd}
\usepackage[mathscr]{eucal}
\usepackage{tikz}

\newcommand{\R}{\mathbb{R}}
\newcommand{\T}{{\bf t}}
\newcommand{\X}{{\bf x}}
\newtheorem{theorem}{Theorem}[section]

\newtheorem{definition}{Definition}[section]

\title{ Hyperplane integrability conditions and smoothing for Radon transforms}

\author{Michael Greenblatt}
\date{\today}

\begin{document}
\maketitle
\begin{abstract} 

This paper may be viewed as a companion paper to [G1]. In that paper, $L^2$ Sobolev estimates derived from a Newton polyhedron-based resolution
of singularities method are combined with interpolation arguments to prove $L^p$ to $L^q_s$ estimates, some sharp up to endpoints, for translation
invariant Radon transforms over hypersurfaces and related operators. Here $q \geq p$ and $s$ can be positive, negative, or zero. 

In this paper, we instead use $L^2$ Sobolev estimates derived from the resolution of singularities methods of [G2]
and combine with analogous interpolation arguments, again resulting in $L^p$ to $L^q_s$ estimates for translation invariant Radon transforms which can be
sharp up to endpoints. It will turn out that
sometimes the results of this paper are stronger, and sometimes the results of [G1] are stronger. 
As in [G1], some of the sharp estimates of this paper occur when $s = 0$, thereby giving new sharp $L^p$ to $L^q$ estimates for 
such operators, again up to endpoints.

Our results lead to natural global analogues whose statements can be recast in terms of a hyperplane integrability condition analogous to that of Iosevich 
and Sawyer in their work [ISa1] on the $L^p$ boundedness of maximal averages over hypersurfaces.

\end{abstract}

\section{ Introduction and theorem statements } 

\subsection{Local $L^p$ to $L^p_s$ boundedness}

We consider the following type of Radon transform operator, defined on functions $\R^{n+1}$, where $\X$ denotes $(x_1,...,x_n)$ and $\T$ denotes $(t_1,...,t_n)$.
\[Tf({\bf x}, x_{n+1}) = \int_{\R^n} f(\X - \T, x_{n+1} - S({\bf t}))\,\phi({\bf t})\,d{\bf t}\tag{1.1}\]
Here $S({\bf t})$ is a real-analytic function on a neighborhood $U$ of the origin, and $\phi({\bf t})$ is a smooth cutoff function supported in $U$. Hence 
$T$ is a convolution operator with a hypersurface measure derived from the graph of $S(\T)$, cutoff by the function $\phi({\bf t})$. By the translation and rotation invariance of the function
space estimates we are proving, without loss of generality we assume that
\[S(0,...,0) = 0 {\hskip 0.65 in} \nabla S(0,...,0) = (0,...,0) \tag{1.2}\]
We also assume $S$ is not identically zero; otherwise $T$ becomes a convolution operator in the first $n$ variables that is easy to analyze.

By the well-known asymptotics for sublevel measures of real analytic functions (we refer to Chapters 6-7 of [AGV] for more details) there are $h> 0$,
an integer $0 \leq d \leq n-1$, and a neighborhood $U$ of the origin in $\R^n$, such that if $V \subset U$ is a neighborhood of the origin then for
some positive constants $b_V, c_V$ one has the following for all $0 < \epsilon < {1 \over 2}$.
\[b_V\epsilon^h |\ln \epsilon|^d < m(\{\T \in V:  |S({\bf t})| < \epsilon\}) < c_V \epsilon^h |\ln\epsilon|^d \tag{1.3}\]
Here $m$ denotes Lebesgue measure. Then our local  $L^p(\R^{n+1})$ to $L^p_{s}(\R^{n+1})$ boundedness theorem is as follows.

\begin{theorem}

There is a neighborhood $W$ of the origin such that if $\phi(\T)$ is supported on $W$ then the following hold.
 
\noindent {\bf 1)} Let $A$ denote the open triangle with vertices $({1 \over 2}, {1 \over n + 1})$, $(0,0)$, and $(1,0)$, and  let $B = \{(x,y) \in A: 
 y < h\}$. Then $T$ is bounded from $L^p(\R^{n+1})$ to $L^p_{s}(\R^{n+1})$ if $({1 \over p}, s) \in B$. 
 
\noindent {\bf 2)} Suppose $h < 1$ and $\phi(\T)$ is nonnegative with $\phi(0) > 0$. Then if $1 < p < \infty$ and 
 $T$ is bounded from $L^p(\R^{n+1})$ to $L^p_{s}(\R^{n+1})$ we must have $s \leq h$.

\end{theorem}

Note that when $h \geq {1 \over n + 1}$, $B$ is just the triangle $A$, while if $h < {1 \over n + 1}$, $B$ is a trapezoid with vertices $(0,0), ({n+1 \over 2}h, h), 
(1 -  {n+1 \over 2}h, h)$, and $(1,0)$. Thus when $h < {1 \over n + 1}$, part b) of Theorem 1.1 shows that part a) gives the sharp amount of $L^p$ to $L^p_s$ improvement up to 
endpoints when $p \in ({n+1 \over 2}h, 1 -  {n+1 \over 2}h)$, while if $h = {1 \over n + 1}$ the same is true for $p = 2$.

Some motivation for Theorem 1.1 is as follows. Let $\rho$ be the measure such that $Tf = f \ast \rho$ in $(1.1)$. Then by the connection between sublevel set measures
and oscillatory integral decay estimates (see Ch 6-7 of [AGuV] for more information), 
the supremum of the numbers $\epsilon$ for which we have a Fourier transform  decay estimate of the form $|\hat{\rho}(0,...,0,\xi)| \leq C(1 + |\xi|)^{-\epsilon}$ for 
all $\phi$ is supported on a sufficiently small neighborhood of the origin is given by the quantity $h$ of $(1.3)$. In any other direction, we get arbitarily fast decay. 

Furthermore,
a straightforward argument shows that the supremum of the $s$ for which one has $L^2(\R^{n+1})$ to $L^2_s(\R^{n+1})$ boundedness of $T$ is exactly the supremum of the
$\epsilon$ for which $|\hat{\rho}(\xi)| \leq C(1 + |\xi|)^{-\epsilon}$ for all $\xi$. Thus the $L^2$ case of Theorem 1.1 says that the estimate for $\hat{\rho}$ that holds in the
$(0,0,...,0,1)$ direction holds uniformly in all directions (up to endpoints), so long as $h \leq {1 \over n +1}$. The interpolation argument we will use
extends this to the whole interval $({n+1 \over 2}h, 1 -  {n+1 \over 2}h)$ when $ h < {1 \over n + 1}$.

\subsection {Global $L^p$ to $L^p_s$ boundedness}

Now suppose $S$ is a compact hypersurface in $\R^{n+1}$ with real analytic boundary. Let $\mu$ denote the standard Euclidean 
surface measure on $S$ and let $U$ be the operator defined on functions on $\R^{n+1}$ by $Uf =  f \ast \mu$. The analogue to Theorem 1.1 for such operators
can be stated in terms of a hypersurface integrability condition analogous to that of Iosevich and Sawyer in their paper [ISa1] on maximal averages. 
Namely, for any hyperplane $P$ in $\R^{n+1}$ we define $\eta_P$ by
\[\eta_P = \sup \{\delta: \int_S (dist(x,P))^{-\delta} d\mu < \infty\} \tag{1.4}\]
Here $dist$ denotes the usual Euclidean distance. We define the index $\eta$ by
\[\eta = \inf_P \eta_P \tag{1.5}\]
Equivalently, $\eta$ is the supremum of the $\delta$ for which  $\int_S (dist(x,P))^{-\delta} d\mu < \infty$ for all hyperplanes $P$. As will be described below,
we must have $\eta > 0$. Our  $L^p(\R^{n+1})$ to $L^p_{s}(\R^{n+1})$ boundedness theorem for $U$ is as follows.

\begin{theorem}
\
\

\noindent {\bf 1)} Let $A$ again denote the open triangle with vertices $({1 \over 2}, {1 \over n + 1})$, $(0,0)$, and $(1,0)$, and  let $D = \{(x,y) \in A: 
 y < \eta\}$. Then $U$ is bounded from $L^p(\R^{n+1})$ to $L^p_{s}(\R^{n+1})$ if $({1 \over p}, s) \in D$.  
 
\noindent {\bf 2)} Suppose $\eta < 1$. If $1 < p < \infty$ and 
 $U$ is bounded from $L^p(\R^{n+1})$ to $L^p_{s}(\R^{n+1})$ then we must have $s \leq \eta$.

\end{theorem}

We refer to Section 1.5 of [G2] for more information concerning the connection between $L^p$ boundedness properties of maximal averages over 
hypersurfaces such as in [ISa1] and the $L^2$ Sobolev smoothing properties of Radon transforms.

Theorem 1.2 follows from Theorem 1.1 in relatively
 short order. For one can use a partition of unity to write $U = \sum_{i=1}^m U_i$, where each $U_i$ is, possibly
after a rotation and translation, of the form $(1.1)$. Furthermore, if  $P$ denotes the tangent plane to $S$ at $x_0$ and $N$ is a sufficiently small neighborhood of the origin, then by the relationship 
between integrals and distribution functions, the index $h$ in Theorem 1.1 corresponding to $x_0$ satisfies
\[ h = \sup \{t: \int_{N \cap S} (dist(x,P))^{-t} d\mu < \infty\}  \tag{1.6}\] 
This is at least as large as 
$\sup \{t: \int_S (dist(x,P))^{-t} d\mu < \infty\} $, which in turn is at least $\eta$. Hence $h \geq \eta$. In other words, for each $i$ the 
analogue of $h$ for the neighborhood $U_i$ is at least  $\eta$.
Therefore adding Theorem 1.1 part 1) over the different $U_i$ gives that  
$U$ is bounded from  $L^p(\R^{n+1})$ to $L^p_{s}(\R^{n+1})$ if $({1 \over p}, s) \in D$ as needed. This gives the first part of Theorem 1.2.

As for the second part of Theorem 1.2, suppose $\eta < 1$ and  suppose $1 < p < \infty$ is such that $U$ is bounded from $L^p(\R^{n+1})$ to 
$L^p_{s}(\R^{n+1})$.
By duality, if $p'$ is such that ${1 \over p } + {1 \over p'} = 1$, we also have that $U$ is bounded 
from $L^{p'}(\R^{n+1})$ to $L^{p'}_{s}(\R^{n+1})$. Interpolating, we get that $U$ is bounded from $L^2(\R^{n+1})$ to $L^2_{s}(\R^{n+1})$.
Thus the measure $\mu$ for which $Uf = f \ast \mu$ satisfies an estimate
\[\hat{\mu}(\xi) \leq C(1 + |\xi|^2)^{-{s \over 2}} \tag{1.7}\]
Let $x_0$ be any point on the boundary of $S$ such that the quantity $h$ of Theorem 1.1 corresponding to $x_0$ satisfies $h < 1$, and let $\psi(x)$ be a nonnegative cutoff function on $\R^{n+1}$ supported near $x_0$ satisfying $\psi(x_0) > 0$. 
Then we must also have $|\widehat{\psi(x)\mu}(\xi)| \leq C_{\psi}(1 + |\xi|^2)^{-{s \over 2}}$.
In other words,  $U_{\psi} f = f \ast \psi(x)\mu$ is bounded from $L^2(\R^{n+1})$ to $L^2_{s}(\R^{n+1})$.
Thus if the support of $\psi$ is sufficiently small, we have that $s$ is at most the quantity $h$ of Theorem 1.1 for the surface $S$ at $x_0$. Taking the infimum over all $x_0$ 
whose corresponding $h$ are less than 1 gives
$s \leq \eta$ as needed. This completes the proof of part 2 of Theorem 1.2.

Note that $\eta$ can never be zero; the proof of the first part of Theorem 1.2 above shows that $U$ has a nonzero amount of 
$L^2$ Sobolev smoothing since $h > 0$ for any $S(\T)$.

\subsection {$L^p$ to $L^q_s$ boundedness for $p \neq q$}

To extend Theorems 1.1 and 1.2 to $L^p$ to $L^q_s$ estimates, we will interpolate with the following consequence of Theorem 1.2 of [G1]. We will use it 
for $p$ approaching $1$ and $q$ approaching infinity. 

\begin{theorem}{\rm ([G1])} For any $1 < p \leq  q < \infty$ and any $\gamma > 1$,
the operator $T$ is bounded from $L^p(\R^{n+1})$ to $L^q_{-\gamma}(\R^{n+1})$.
\end{theorem}

Let $Q$ be the plane in $\R^3$ containing the line $\{(x,y,z): x = y, z = h\}$ and the point 
$(1,0, -1)$, and let $Q'$ be the plane containing the line $\{(x,y,z): x = y, z = \eta\}$ and the point $(1,0, -1)$.

Suppose $h < {1 \over n + 1}$ and we are in the setting of Theorem 1.1. If we interpolate the estimate from Theorem 1.1 corresponding to the
upper edge of the trapezoid $B$ with Theorem 1.2, and then let $(p,q)$ go to $(1,\infty)$, we obtain an $L^p$ to $L^q_s$ 
boundedness theorem for $({1 \over p}, {1 \over q}, s)$ below a triangle $Y$  in the plane $Q$. The analogous statement holds if $\eta <
 {1 \over n + 1}$ and we are in the setting of Theorem 1.2, where $Q$ is replaced by $Q'$. We can then interpolate these  results
 with the trivial $L^p$ to $L^p$ boundedness of Radon transforms for $1 < p < \infty$ to obtain the following analogue of the first part 
of Theorem 1.4 of [G1], keeping in mind that if $s_1 < s_2$ then $L^q_{s_2}(\R^{n+1}) \subset L^q_{s_1}(\R^{n+1})$ continuously for any 
$1 < q < \infty$.  

\begin{theorem} 

\
\
 
Suppose we are in the setting of Theorem 1.1 with $h < {1 \over n + 1}$. Let $Y$ be the closed triangle with vertices $({n+1 \over 2}h,
 {n+1 \over 2}h, h)$, $(1 - {n+1 \over 2}h, 1 - {n+1 \over 2}h, h)$, and $(1,0,-1)$.
Let $Y_1$ be the closed triangle with vertices $(0,0,0)$, $(1,0, -1)$, and $({n+1 \over 2}h, {n+1 \over 2}h, h)$ and  let $Y_2$ be the closed triangle
with vertices $(1,1,0)$, $(1,0,-1)$, and $(1 - {n+1 \over 2}h, 1 - {n+1 \over 2}h, h)$. 

There is a neighborhood $W$ of the origin such that if $\phi(\T)$ is supported on $W$, then if $({1 \over p}, {1 \over q}, s)$ is such that there is a $t > s$ with $({1 \over p}, {1 \over q}, t)$ is in the interior of $Y \cup Y_1 \cup Y_2$, then $T$ is bounded from $L^p(\R^{n+1})$ to $L^q_s(\R^{n+1})$. 

The analogous statement holds in the setting Theorem 1.2 if $T$ is replaced by $U$ and $h$ is replaced by $\eta$.
 
\end{theorem}

The triangles $Y$, $Y_1$, and $Y_2$ in Theorem 1.4 can be visualized as follows. The segment from $({n+1 \over 2}h, {n+1 \over 2}h, h)$ to 
$(1 - {n+1 \over 2}h, 1 - {n+1 \over 2}h, h)$ is a line segment above the line $y = x$, at fixed height $z = h$, which is symmetric
about the midpoint $({1 \over 2}, {1 \over 2}, h)$. The trangle $Y$ is then the convex hull of this segment with the point $(1,0,-1)$ that is below the lower-rightmost 
point in the square $[0,1] \times [0,1]$.
 $Y_1$ is the convex hull of the left side of $Y$ with $(0,0,0)$ and $Y_2$ is  the convex hull of the right side of $Y$ with $(1,1,0)$, so 
that $Y_1$ and $Y_2$ are symmetric about the plane $x + y = 1$. 

If $h \geq {1 \over n+1}$ or $\eta \geq {1 \over n+1}$ in Theorem 1.1 or 1.2 respectively, the analogous interpolation of
Theorem 1.1 or 1.2 with Theorem 1.3 and the trivial $L^p$ to $L^p$ statements gives the following analogue of the second part of Theorem 1.4 of [G1].

\begin{theorem}
Suppose we are in the setting of Theorem 1.1 and $h \geq {1 \over n + 1}$. Let $Y_3$ be the closed triangle with vertices $(0,0,0)$, $(1,0,-1)$, and $({1 \over 2}, {1 \over 2}, {1 \over n + 1})$ and let $Y_4$  be the closed triangle with vertices $(1,1,0)$, $(1,0,-1)$, and $({1 \over 2}, {1 \over 2}, {1 \over n + 1})$.
There is a neighborhood $W$ of the origin such that if $\phi(\T)$ is supported on $W$ then if $({1 \over p}, {1 \over q}, s)$ is such that there is a $t > s$ with $({1 \over p}, {1 \over q}, t)$ is in the interior of $Y_3 \cup Y_4$, then $T$ is bounded from $L^p(\R^{n+1})$ to $L^q_s(\R^{n+1})$. 

If we are in the setting of Theorem 1.2 and $\eta \geq {1 \over n + 1}$, the analogous statement holds with $T$ replaced by $U$ and the condition $h \geq
{1 \over n + 1}$ replaced by $\eta \geq {1 \over n + 1}$.

\end{theorem}

So in Theorem 1.5, $Y_3$ is the convex hull of $(0,0,0)$ and the line segment from $({1 \over 2}, {1 \over 2}, {1 \over n+1})$ to $(1,0,-1)$ on the plane $x + y = 1$, and $Y_4$ is the convex hull of $(1,1,0)$ and the line segment from 
$({1 \over 2}, {1 \over 2}, {1 \over n+1})$ to $(1,0,-1)$. Again, $Y_3$ and $Y_4$ are symmetric about the plane 
$x + y = 1$, this time with a common edge on this plane.

\subsection{Sharpness when $p \neq q$}

Theorem 1.5 will never be sharp up to endpoints when $h > {1 \over n + 1}$ since the portion of Theorem 1.1 being used is not sharp. However, 
Theorem 1.4 is sharp in a number of situations. To help understand when, we need some terminology and results from [G1].

\begin{definition}

Let $f({\bf t})$ be a real analytic function defined on a neighborhood of the origin in 
$\R^n$, and
let $\sum_{\alpha} f_{\alpha}{\bf t}^{\alpha}$ denote the Taylor expansion of $f({\bf t})$ at the origin.
For any $\alpha$ for which $f_{\alpha} \neq 0$, let $Q_{\alpha}$ be the octant $\{{\bf t} \in \R^n: 
t_i \geq \alpha_i$ for all $i\}$. Then the {\it Newton polyhedron} $N(f)$ of $f({\bf t})$ is defined to be 
the convex hull of all $Q_{\alpha}$.  

\end{definition}

\begin{definition} 

Where $f({\bf t})$ is as in Definition 1.1, define $f^*({\bf t})$ by 
$$f^*({\bf t}) = \sum_{(v_1,...,v_n)\,\,a \,\,vertex \,\,of\,\,N(f)} |t_1|^{v_1}...|t_n|^{v_n} \eqno (1.8)$$
\end{definition}

\begin{definition}Suppose $F$ is a compact face of the Newton polyhedron $N(f)$. Then
if $\sum_{\alpha} f_{\alpha}{\bf t}^{\alpha}$ denotes the Taylor expansion of $f$ like above, 
define $f_F({\bf t}) = \sum_{\alpha \in F} f_{\alpha}{\bf t}^{\alpha}$.

\end{definition}

\begin{definition} For $f({\bf t})$ as above, we denote by $o(f)$ the maximum order of any zero of any $f_F({\bf t})$ on $(\R - \{0\})^n$. We take
$o(f) = 0$ if there are no such zeroes.

\end{definition}

\begin{definition} The Newton distance $d(f)$ is defined to be the minimal $t$ for which $(t,...,t)$ is in the Newton
polyhedron $N(f)$.

\end{definition}

By Lemma 2.1 of [G3], similarly to $(1.3)$ there is an $r_0 > 0$, an $g> 0$, and an integer $d_0$ satisfying $0 \leq d_0 \leq n-1$, such that if $r < r_0$ then there are positive constants $b_r$ and $B_r$ such that for $0 < \epsilon < {1 \over 2}$ we have
\[b_r \epsilon^g |\ln \epsilon|^{d_0} < m(\{{\bf t} \in (0,r)^n:  S^*({\bf t}) < \epsilon\}) < B_r \epsilon^g |\ln\epsilon|^{d_0}\tag{1.9}\]
Here $m$ denotes Lebesgue measure. The sharp portion of Theorem 1.4 of [G1] in the setting of Theorem 1.1 is as follows.

\begin{theorem}{\rm ([G1])} There is a neighborhood
 $W$ of the origin such that if $\phi(\T)$ is supported on $W$ then the following hold.
 
Suppose $g < {1 \over \max(o(S),2)}$, where $o(S)$ is as in Definition 1.4. Let $Z$ be the closed triangle with vertices $({\max(o(S), 2) \over 2}g, {\max(o(S), 2) \over 2}g, g)$,
$(1 - {\max(o(S), 2) \over 2}g, 1 - {\max(o(S), 2) \over 2}g, g)$, and $(1,0,-1)$. Let $Z_1$ be the closed triangle with vertices $(0,0,0)$, $(1,0, - 1)$,
 and $({\max(o(S), 2) \over 2}g, {\max(o(S), 2) \over 2}g, g)$, and let $Z_2$ be the closed triangle
with vertices $(1,1,0)$, $(1,0,-1)$, and $(1 - {\max(o(S), 2) \over 2}g, 1 - {\max(o(S), 2) \over 2}g, g)$. If $({1 \over p}, {1 \over q}, s)$ is such that there is a $t > s$ with $({1 \over p}, {1 \over q}, t)$ in the interior of $Z \cup Z_1 \cup Z_2$, then $T$ is bounded from $L^p(\R^{n+1})$ to $L^q_s(\R^{n+1})$. 
\end{theorem}

Let $P$ be the plane with equation $(g + 1)(x - y) + z = g$. Then the triangle $Z$ is a subset of $P$. Theorem 1.6 is sharp in the following sense, 
as shown in [G1].

\begin{theorem}{\rm ([G1])} Suppose $g \leq {1 \over \max(o(S),2)}$ and there is a $C_1 > 0$ and a 
neighborhood $N_0$ of the origin such that $\phi(\T) > C_1$ on $N_0$. Then for any $1 < p, q < \infty$, 
if $({1 \over p}, {1 \over q}, s)$ is 
such that there is a $t < s$ with $({1 \over p}, {1 \over q}, t)$ on the plane $P$, then $T$ is not bounded from $L^p(\R^{n+1})$ to $L^q_s(\R^{n+1})$.
\end{theorem}
 
\noindent To help us understand the relation between the results of [G1] and of this paper, we first observe that by Lemma 2.1 of [G4], for any
 real analytic function $f({\bf t})$ defined near the origin there is a neighborhood
of the origin on which $|f(\T)| \leq Cf^*(\T)$ for some $C > 0$. Hence by the definitions $(1.3)$ and $(1.9)$ of $g$ and $h$ one always has $h \leq g$.
As a result, if $n+1 \geq \max(o(S),2)$, the results of this paper are contained in that of [G1] as the vertices of the various triangles in Theorem 1.4
will be no higher than the corresponding triangles in Theorem 1.6. But when $n+1 < \max(o(S),2)$ this does not have to be the case. 

Suppose $n + 1 < \max(o(S),2)$. It follows from Theorem 1.2 of [G4] that if $o(S) \leq  d(S)$, where $d(S)$ and $o(S)$ are as defined above, then one has
$g = {1 \over d(S)} = h$. Comparing the statements of Theorems 1.4 and 1.6 when  $n + 1 < \max(o(S), 2)$ and $g = h$, we see that
the plane $P$ containing the triangle $Z$ in Theorem 1.6 is the same as the plane $Q$ containing the triangle $Y$ in Theorem 1.4, and that 
the statement of Theorem 1.4 is strictly stronger than that of Theorem 1.6 whenever the assumptions of both are satisfied; this is because $Z \subsetneq Y$. 

If in the $n + 1 < \max(o(S),2)$, $o(S) \leq  d(S)$ scenario we wish to use the sharpness result, Theorem 1.7, in conjunction with Theorem 1.4, then we must also have that $g \leq {1 \over \max(o(S),2)}$ since that is assumed in Theorem 1.7. Once this assumption is added, Theorem 1.7 says that for $({1 \over p}, {1 \over q}, s)$ above $P = Q$ one does not have $L^p$ to $L^q_s$ boundedness. Hence for $({1 \over p}, {1 \over q}, s)$
below the interior of $Y$ we have boundedness, and for $({1 \over p}, {1 \over q}, s)$ above the interior of $Y$ we do not. 
Thus the level of Sobolev improvement is sharp 
up to endpoints for $({1 \over p}, {1 \over q})$ in the projection of the interior of $Y$ onto the $x$-$y$ plane, under our assumptions that $o(S) \leq  d(S)$, 
$n+1 < \max(o(S),2)$, and $g \leq {1 \over \max(o(S),2)}$. Examples were these assumptions are satisfied are easy to construct.

If $n+1 < \max(o(S),2)$ but $o(S) > d(S)$, then it is still possible that $g = h$, in which case considerations similar to the above again apply, but in the 
more common situation that $g > h$, we do not have such statements.

It is worth pointing out that Theorem 1.4 of [G1] also has a non-sharp portion analogous to Theorem 1.5 of this paper. Under the hypotheses of Theorem 1.5,
the statement of this part of Theorem 1.4 of [G1] gives the statement analogous to that of Theorem 1.5 where the condition $h \geq {1 \over n + 1}$ is
replaced by the condition $g \geq{1 \over \max(o(S),2)}$ and the upper vertex 
$({1 \over 2},  {1 \over 2}, {1 \over n + 1})$ of $Y_3$ and $Y_4$ is replaced by $({1 \over 2}, {1 \over 2}, {1 \over \max(o(S),2)})$. Thus once again,
if $h = g$ we have a stronger result in this paper when $n + 1 < \max(o(S),2)$, and a stronger result in [G1] when $n + 1 > \max(o(S),2)$.

\section{Some background}

There has been quite a bit of work done on the boundedness properties of Radon transforms on function spaces, so we focus our attention on 
Sobolev space improvement and $L^p$ to $L^q$ improvement results for Radon transforms over hypersurfaces. For curves in $\R^2$, [S] provides comprehensive $L^p_{\alpha}$ to 
$L^q_{\beta}$ boundedness results for Radon transforms that are sharp up to endpoints. These results include general non-translation invariant operators.

For translation invariant Radon transforms, $L^2$ to $L^2_{\beta}$ Sobolev space improvement is equivalent to a surface measure Fourier transform 
decay rate estimate. When $n = 2$, the stability theorems of Karpushkin [Ka1] [Ka2] combined with [V] give such sharp decay rate results for real
analytic surfaces. Extensions to finite type smooth surfaces appear in [IkKM]. For general $n$, the author's paper [G2] provides such some such estimates,
and there are also some earlier results for example for convex (non necessarily smooth) $S(\T)$ such as [R1][R2].  For general $p$, the paper [St]  considers Sobolev estimates for Radon transforms in a quite general setting, focusing attention on singular $\phi({\T})$.

For $L^p$ to $L^q$ improvement for Radon transforms over hypersurfaces, there have been a number of other results for 
Radon transforms. The case of surfaces with nonvanishing Gaussian curvature are covered in [L][Ste][Str]. 
The situation where the  $S({\bf t})$ is a homogeneous or mixed homogeneous function has been considered in [DZ] [FGU1] [FGU2] [ISa2].
Convex surfaces were considered in [ISaS]. Also, there have been papers considering damped Radon transforms, where instead of the function
$\phi(\T)$ in $(1.1)$ one uses $\phi(\T)\psi(\T)$, where $\psi(\T)$ has zeroes on a set chosen to be natural for the surfaces at hand. Often $\psi(\T)$
is related to the Hessian determinant of $S(\T)$. We mention [Gr] and [O] as examples of such results.

\section{Examples}

\noindent {\bf Example 1.} 

We consider the case where $n = 1$, so that we are considering Radon transforms over curves in the plane. In this situation we have
$S(t) = ct^l + O(t^{l+1})$ for some integer $l \geq 2$, where $c \neq 0$. The index $h$ is determined by $(1.3)$ to be $h = {1 \over l}$. Hence $h < {1 \over n+1} = {1 \over 2}$ except when
$l = 2$, in which case $h = {1 \over n + 1}$. In the former case, $B$ is the trapezoid with vertices $(0,0), ({1 \over l}, {1 \over l}), ({l - 1 \over l}, {l - 1 \over l})$, and $(1,0)$,
and in the $l = 2$ case $B$ is just the triangle $A$, whose vertices are $(0,0), ({1 \over 2}, {1 \over 2}),$ and $(1,1)$. This is the same range of $L^p$ to $L^p_s$ boundedness following from [G3]. This also follows from [C], where it is also shown that if $l > 2$ one has $L^p$ to $L^p_s$ boundedness
on the boundary of $B$ except at the vertices $({1 \over l}, {1 \over l})$ and $({l - 1 \over l}, {l - 1 \over l})$, and that one does not have $L^p$ to 
$L^p_s$ boundedness for $s > 0$ outside the closure of $B$.

Moving to what Theorems 1.4 and 1.5 say here, Theorem 1.4 applies when $l > 2$ and Theorem 1.5 applies when $l = 2$. If $l > 2$ the triangle $Y$ has vertices $({1 \over l},
{1  \over l}, {1 \over l})$, $({l - 1 \over l},  {l - 1\over l}, {1 \over l})$, and $(1,0, - 1)$. The triangle $Y_1$ has vertices $(0,0,0)$, $(1,0,-1)$, and 
$({1  \over l},{1 \over l},{1 \over l})$, and $Y_2$ has vertices $(1,1,0)$, 
$(1,0, -1 )$, and $({l - 1 \over l}, {l - 1 \over l}, {1 \over l})$. Theorem 1.4 then says that one
has $L^p$ to $L^q_s$ boundedness for $({1 \over p}, {1 \over q}, s)$ below the interior of $Y \cup Y_1 \cup Y_2$, and Theorem 1.7 says one does not have $L^p$ to $L^q_s$ boundedness for $({1 \over p}, {1 \over q}, s)$ above the plane $P$ containing $Y$.

In the case where  $l = 2$ the triangle $Z$ reduces to a line and Theorem 1.5 applies. This time $Y_3$ has vertices $(0,0,0), (1,0,-1)$, and $({1 \over 2}, {1 \over 2}, {1 \over 2})$
and $Y_4$ has vertices $(1,1,0)$, $(1,0,-1)$, and $({1 \over 2}, {1 \over 2}, {1 \over 2})$. Theorem 1.5 now says that one
has $L^p$ to $L^q_s$ boundedness for $({1 \over p}, {1 \over q}, s)$ below the interior of $Z_3 \cup Z_4$, and Theorem 1.7 says one does not 
have $L^p$ to $L^q_s$ boundedness for 
$({1 \over p}, {1 \over q}, s)$ above the plane $P$, which intersects $Y_3 \cup Y_4$ in the line segment connecting $({1 \over 2}, {1 \over 2}, {1 \over 2})$ to $(1,0,-1)$.

The above results are exactly the same results given for curves in the plane provided by Theorem 1.6 from [G1] since $h = g = {1 \over l}$ and $\max(o(S),2) = 2 = n + 1$ here.

\noindent {\bf Example 2.}

Suppose now that $n > 1$. We write out the $L^p$ to $L^q$ estimates that follow from Theorem 1.4. Assume that $h < {1 \over n + 1}$ so that 
the hypotheses of Theorem 1.4 hold. We focus our attention on the intersection of the triangle 
$Y$ with the $x$-$y$ plane. The line segment with vertices $({n+1 \over 2}h, {n+1 \over 2}h, h)$ and $(1,0,-1)$ intersects the $x$-$y$ plane at
\[{1 \over h + 1}\bigg({n+1 \over 2}h, {n+1 \over 2}h, h\bigg) + {h \over h + 1} (1,0,-1)\]
\[ = \bigg( {h (n + 3) \over 2(h+1)}, {h (n + 1) \over 2(h+1)}, 0\bigg) \tag{3.1}\]
The line segment with vertices $(1 - {n+1 \over 2}h, 1 - {n+1 \over 2}h, h)$ and $(1,0,-1)$ intersects the $x$-$y$ plane at 
\[{1 \over h + 1}\bigg(1 - {n+1 \over 2}h, 1 - {n+1 \over 2}h, h\bigg) + {h \over h + 1} (1,0,-1) \]
\[=\bigg(1 -{h (n + 1) \over 2(h+1)}, 1 -{h (n + 3) \over 2(h+1)}, 0\bigg)\tag{3.2}\]
Thus if we let $p_1 = ( {h (n + 3) \over 2(h+1)}, {h (n + 1) \over 2(h+1)}, 0)$ and $p_2 = (1 -{h (n + 1) \over 2(h+1)}, 1 -{h (n + 3) \over 2(h+1)}, 0)$,
the intersection of the triangle $Y$ with the $x$-$y$ plane is the line segment from $p_1$ to $p_2$. One can check similarly to the above that the intersection of the triangle $Y_1$ of Theorem 1.4 with the $x$-$y$ plane is the line segment from $(0,0,0)$ to $p_1$, and the intersection of the triangle $Y_2$
 with the $x$-$y$ plane is the line segment from $p_2$ to $(1,1,0)$. 
 
 Let $J$ be the trapezoid in two dimensions with vertices $(0,0),  ( {h (n + 3) \over 2(h+1)}, {h (n + 1) \over 2(h+1)}), (1 -{h (n + 1) \over 2(h+1)}, 1 -{h (n + 3) \over 2(h+1)})$, and $(1,1)$. Then Theorem 1.4 gives $L^p$ to $L^q$ boundedness for $T$ when $({1 \over p}, {1 \over q})$ is in the interior
 of $J$. 
 
Next, we add the assumption that $o(S) \leq d(S)$, where $d(S)$ is the Newton distance of $S$ as in Definition 1.5. Then as indicated in the
discussion at the end of section 
1, Theorem 1.2 of [G4] implies that $g = h = {1 \over d(S)}$. So by that discussion, if $\max(o(S),2) \leq n + 1$, Theorem 1.6 gives a result at least as
 strong as Theorem 1.4, while if $\max(o(S),2) > n + 1$ then Theorem 1.4 gives a strictly stronger result than Theorem 1.6. Furthermore, substituting $h = {1 \over d(S)}$,  the trapezoid $J$ can be described in terms of $d(S)$ instead of $h$ as the trapezoid with vertices
 $(0,0),  ( {n + 3 \over 2d(S) + 2}, {n + 1 \over 2d(S) + 2}), (1 -{n + 1 \over 2d(S) + 2}, 1 -{n + 3 \over 2d(S) + 2})$, and $(1,1)$.
 
 If we add the additional assumption
 that $g \leq {1 \over \max(o(S),2)}$, then Theorem 1.7 now applies and one cannot have $L^p$ to $L^p_s$ boundedness for $({1 \over p}, {1 \over q},
 s)$ above the plane $P$ of that theorem, which is now the same as the plane $Q$ containing the triangle $Y$ in Theorem 1.4 since $g = h$. In terms of $L^p$ to $L^q$ estimates, this means that one cannot have $L^p$ to $L^q$ estimates for $({1 \over p}, {1 \over q})$ below the line containing the segment from 
$ ( {n + 3 \over 2d(S) + 2}, {n + 1 \over 2d(S) + 2})$ to $(1 -{n + 1 \over 2d(S) + 2}, 1 -{n + 3 \over 2d(S) + 2})$, which has equation $y = x - {1 \over 
d(S) + 1}$.
 
\section{The proof of Theorem 1.1}

We consider the operator $T$ in the form $Tf = f \ast \rho$, where $\rho$ is a hypersurface measure supported near the origin.
 In Section 2 of [G2], the Fourier transform $\hat{\rho}(\xi)$ was analyzed. It is shown there that there is a neighborhood of
$W$ of the origin such that if $\phi(\T)$ in $(1.1)$ is supported in $W$, then one can do the following. Up to a set of measure zero, one 
can write (independently of $\phi$) $\R^n = \cup_{i=1}^N K_i$, where each $K_i$ is a wedge in $\xi$ space of the form 
$K_i = \{\xi \in \R^n: a_{ik} \cdot \xi > 0$ for $k = 1,...,M_i\}$ such that $K_i$ have the following properties.

For a given $i$, one can write $\phi(\T) = \sum_{j=1}^{N_i}\phi_{ij}(\T)$, where $\sup_{\T} |\phi_{ij}(\T)| \leq sup_{\T}|\phi(\T)|$, such that we have 
the following. Let $T_{ij}$ be the operator
\[T_{ij}f({\bf x}, x_{n+1}) = \int_{\R^n} f(\X - \T, x_{n+1} - S({\bf t}))\phi_{ij}({\bf t})\,d{\bf t}\tag{4.1}\]
Thus $\sum_{j=1}^{N_i} T_{ij} = T$.
Let the measure $\rho_{ij}$ be such that $T_{ij}f = f \ast \rho_{ij}$, so that $\sum_{j=1}^{N_i} \rho_{ij} = \rho$.
By $(2.25)$ of [G2],  for $\xi \in K_i$,  if $r$ denotes the diameter of $W$ one has an estimate 
\[|\widehat{\rho_{ij}}(\xi)| \leq C\int_{\{x: |x| < r\}}\min(1, (|\xi||g_{ij}(x)|)^{-{1 \over n+1}})\,dx \tag{4.2}\]
Here $g_{ij}$ is either $S(x)$ or a linear function that is zero at the origin. Furthermore, the subdivisions of section 2 of [G2] are such that each
$\phi_{ij}(\T)$ can be written as an infinite sum $\sum_{k=0}^{\infty} \phi_{ijk}(\T)$, where $\phi_{ijk}(\T)$ is nonzero only if $c_12^{-k} < |g_{ij}(\T)| < c_22^{-k}$ for some
$c_2 > c_1 > 0$, such that
one has the following analogue of $(4.2)$. Define $T_{ijk}$ by
\[T_{ijk}f({\bf x}, x_{n+1}) = \int_{\R^n} f(\X - \T, x_{n+1} - S({\bf t}))\phi_{ijk}({\bf t})\,d{\bf t}\tag{4.3}\]
Let $\rho_{ijk}$ be such that $T_{ijk}f = f \ast \rho_{ijk}$, so that $\sum_k \rho_{ijk} = \rho_{ij}$. Then the arguments leading to $(2.25)$ in [G2] are such that $(4.2)$ is obtained by adding the analogous estimates for $\rho_{ijk}$, namely that whenever $\xi \in K_i$ one has
\[|\widehat{\rho_{ijk}}(\xi)| \leq C\int_{\{x: |x| < r,\,\,c_12^{-k} <  |g_{ij}(x)| < c_2 2^{-k} \}}\min(1, (|\xi||g_{ij}(x)|)^{-{1 \over n+1}})\,dx \tag{4.4}\]
In view of the above, it is natural to write $T = \sum_{ijk} T_{ijk}$ as follows. Let $M_if$ be the Fourier multiplier operator with multiplier given by
$\chi_{K_i}(\xi)$, and then let $U_{ijk}f = M_i f \ast \rho_{ijk}$. Note that  $\sum_{jk}U_{ijk}f = M_if \ast \sum_{jk}\rho_{ijk} = M_i Tf$, so that we have
\[\sum_{ijk}U_{ijk}f = \sum_i M_i T f = Tf\tag{4.5}\]
Note that if we write $U_{ijk}f = f \ast \sigma_{ijk}$ for a distribution $\sigma_{ijk}$, then $(4.4)$ combined
with the support condition in $\xi$ space induced by $\chi_{K_i}(\xi)$ give that for {\it any} $\xi$ one has 
\[|\widehat{\sigma_{ijk}}(\xi)| \leq C\int_{\{x: |x| < r,\,\,c_12^{-k} <  |g_{ij}(x)| < c_2 2^{-k} \}}\min(1, (|\xi||g_{ij}(x)|)^{-{1 \over n+1}})\,dx \tag{4.6}\]
We now move to the proof of Theorem 1.1. We break into cases $h \leq {1 \over n +1}$ and $h > {1 \over n+1}$, starting with the former case.

Suppose that $h \leq {1 \over n+1}$. We examine $(4.6)$ in this situation. The main issue is the case where $g_{ij}(x)$ is 
$S(x)$. So we assume $g_{ij}(x) = S(x)$ for now and we examine what happens when $(1.3)$ is inserted into $(4.6)$. We obtain that for any $h' < h$
and any $\xi$ we have
\[|\widehat{\sigma_{ijk}}(\xi)| \leq C'\min(1, |\xi|^{-{1 \over n+1}} 2^{{k \over n+1}})\times 2^{-kh'}\tag{4.7}\]
\[\leq C'\min(2^{k \over n+1}, |\xi|^{-{1 \over n+1}} 2^{{k \over n+1}})\times 2^{-kh'}\]
\[\leq C'' (1 + |\xi|)^{-{1 \over n+1}}\times 2^{{k \over n+1}}\times 2^{-kh'}\tag{4.8}\]
Hence we have the bound
\[||U_{ijk}f||_{L^2_{1 \over n+1}} \leq C'''\times 2^{k \over n+1}\times 2^{-kh'}||f||_{L^2} \tag{4.9}\]
On the other hand, $U_{ijk}f = M_i f \ast \rho_{ijk}$ is the convolution of $M_i f$ with the measure $f$ is being convolved with in $(4.3)$.
Since the function $\phi_{ijk}(\T)$ is supported in the set where  $c_12^{-k} < |S(\T)| < c_22^{-k}$,  this is the convolution of
$M_i f$ with a finite measure of total size bounded by $C2^{-kh'}$ for any $h' < h$. 
Hence for any $1 < p < \infty$, $||U_{ijk}f||_p \leq C2^{-kh'}||M_if||_p$. However,
since $M_i$ is a multiplier operator whose multiplier is given by the characteristic function of a wedge defined by hyperplaces, we also have 
$||M_if||_p \leq C'||f||_p$. We conclude that for any $i,j,$ and $k$ we have
\[||U_{ijk}f||_p \leq C''''2^{-kh'}||f||_p \tag{4.10}\]
Notice that since $h' < h \leq  {1 \over n+1}$, the power of $2$ in $(4.9)$ is increasing in $k$, while the power of $2$ in $(4.10)$ decreases in $k$.
The idea now is to interpolate $(4.9)$ and $(4.10)$ for $p \rightarrow \infty$, to get the strongest result possible that still has a $2^{-\epsilon k}$ 
coefficient for some small $\epsilon > 0$, so that we may sum in $k$.
We observe that the equation $\alpha({1 \over n + 1} - h') + (1 - \alpha)(-h') = 0$ is solved by $\alpha = h'(n+1)$, so that the borderline case where 
$\epsilon = 0$ occurs with this weighting in the interpolation, if we set $p = \infty$. If we let ${1 \over p^* } = \alpha{1 \over 2} + (1 - \alpha)*0 
= {h'(n+1) \over 2}$, and $s^* = \alpha{1 \over n + 1}  + (1-\alpha)*0 = h'$, we obtain the borderline $({1 \over p^*}, s^*) = ({h'(n+1) \over 2}, h')$.

We now take $h' = h_m$ for a sequence
$h_m$ be a sequence increasing to $h$. We let $\alpha_m$ correspondingly converge to $h(n + 1)$ and $p_m$ go to $\infty$ such that 
the power of 2 in the above interpolation decreases in $k$ for $h' = h_m$, $\alpha = \alpha_m$, and $p = p_m$. Thus for each $m$
we have the following estimate, where $s_m = {\alpha_m \over n + 1}$ and where $\epsilon_m > 0$.
\[||U_{ijk}f||_{L^{p_m}_{s_m}} \leq C2^{-\epsilon_m k}||f||_{p_m} \tag{4.11}\]
Define $U_{ij}f = M_if \ast \rho_{ij}$, so that $\sum_k U_{ijk} = U_{ij}$ and we have the finite sum $\sum_{ij}U_{ij} = \sum_iM_i =  T$.
Adding $(4.11)$ over all $k$ then gives $||U_{ij} f||_{L^{p_m}_{s_m}} \leq C'||f||_{p_m}$.
Note that as $m$ goes to infinity, $({1 \over p_m}, s_m)$ converges to the upper left endpoint of the trapezoid $B$ in part 1 of Theorem 1.1, except at 
the borderline case $h = {1 \over n + 1}$ where it converges to the upper vertex of $B$ which is now a triangle. 

If $h < {1 \over n + 1}$ and we
 do the analogous procedure to the above, letting $p \rightarrow 1$ instead of $p \rightarrow \infty$, we 
obtain the analogous statement for the upper right endpoint of this trapezoid. Alternatively, one can just use duality to go from the upper left endpoint 
to the upper right endpoint. 

At any rate, when $h < {1 \over n + 1}$, for a given $m$ we interpolate the estimate $(4.11)$ with the corresponding estimate for the upper right endpoint of the trapezoid, and then
interpolate the result with the trivial $L^p$ to $L^p$ estimates. Letting $m$ go to infinity for both endpoints, we obtain $L^p$ to $L^p_s$ boundedness for $({1 \over p}, s) \in B$ for $U_{ij}$ when $h < {1 \over n + 1}$. If $h = {1 \over n + 1}$, then we do the same procedure 
with just the upper vertex of $B$, now a triangle, and once again get $L^p$ to $L^p_s$ boundedness for $({1 \over p}, s) \in B$.

The above argument was for the (main) case when $g_{ij}(x)= S(x)$. The other possibility is that $g_{ij}(x)$ is a linear function with  a zero at the origin,
where
the index analogous to $h$ for $g_{ij}(x)$ in this case is $1 > {1 \over n+1}$. So the argument of the $h = {1 \over n + 1}$ situation now applies;
the estimates used above from $(1.3)$ from the $h = {1 \over n + 1}$ case will still hold and the argument above goes through.
So once again we have $L^p$ to $L^p_s$ boundedness for $({1 \over p}, s) \in B$ for $U_{ij}$ since the $B$ for the $h = {1 \over n + 1}$ case 
contains $B$ for all $h < {1 \over n + 1}$ cases.

Since $T$ is the (finite) sum $T = \sum_{ij} U_{ij}$, we conclude from the above that $T$ is bounded from  $L^p$ to $L^p_s$  for all  $({1 \over p}, s) \in B$. This completes the proof of part 1 of Theorem 1.1 when $h \leq {1 \over n + 1}$.

Suppose now $h > {1 \over n + 1}$. As in the case above where $g_{ij}(x)$ is a linear function, each $g_{ij}(x)$ will satisfy estimates from $(1.3)$ that 
are stronger than that of the $h = {1 \over n + 1}$ situation. Thus like above each $U_{ij}$ is bounded from $L^p$ to $L^p_s$ for $({1 \over p}, s)$ in the interior of the triangle corresponding to the $h = {1 \over n + 1}$ situation, which is exactly $B$ here. Adding over the finitely many $i$ and $j$ we once 
again have that $T$ is bounded from  $L^p$ to $L^p_s$  for all  $({1 \over p}, s) \in B$. This completes the proof of part 1 of Theorem 1.1.

As for part 2 of Theorem 1.1, suppose $g < 1$ and $\phi(\T)$ is nonnegative with $\phi(0) > 0$. Suppose  further that $1 < p < \infty$ is 
such that $T$ is bounded 
from $L^p(\R^{n+1})$ to $L^p_s(\R^{n+1})$. Then by duality, $T$ is bounded  from $L^{p'}(\R^{n+1})$ to $L^{p'}_s(\R^{n+1})$, where ${1 \over p}
+{1 \over p'} = 1$. Interpolating, we get that $T$ is bounded from $L^2(\R^{n+1})$ to $L^2_s(\R^{n+1})$. In order for part 2) of Theorem 1.1 to hold 
we must verify that $s \leq h$. This however is an immediate consequence of part 4 of Theorem 1.2 of [G2]. This concludes the proof of Part 2 of Theorem 1.1, and therefore the proof of the whole theorem.

\section{References.}

\noindent [AGuV] V. Arnold, S. Gusein-Zade, A. Varchenko, {\it Singularities of differentiable maps},
Volume II, Birkhauser, Basel, 1988. \parskip = 4pt\baselineskip = 3pt

\noindent [C] M. Christ, {\it Failure of an endpoint estimate for integrals along curves} in Fourier analysis and partial differential equations (Miraflores de la Sierra, 1992), 163-168, Stud. Adv. Math., CRC, Boca Raton, FL, 1995. 

\noindent [DZ] S. Dendrinos, E. Zimmermann, {\it On $L^p$-improving for averages associated to mixed homogeneous polynomial hypersurfaces in $\R^3$},
 J. Anal. Math. {\bf 138} (2019), no. 2, 563-595. 

\noindent [FGoU1] E. Ferreyra, T. Godoy, M. Urciuolo, {\it Boundedness properties of some convolution operators with singular measures}, Math. Z. 
{\bf 225} (1997), no. 4, 611-624.

\noindent [FGoU2] E. Ferreyra, T. Godoy, M. Urciuolo, {\it Sharp $L^p$-$L^q$ estimates for singular fractional integral operators.}
Math. Scand. {\bf 84} (1999), no. 2, 213-230. 

\noindent [G1] M. Greenblatt, {\it Smoothing theorems for Radon transforms over hypersurfaces and related operators}, submitted.

\noindent [G2] M. Greenblatt, {\it Fourier transforms of indicator functions, lattice point discrepancy, and the stability of integrals}, submitted.

\noindent [G3] M. Greenblatt, {\it $L^p$ Sobolev regularity of averaging operators over hypersurfaces and the Newton polyhedron}, J. Funct. Anal. 
{\bf 276} (2019), no. 5, 1510-1527.

\noindent [G4] M. Greenblatt, {\it Oscillatory integral decay, sublevel set growth, and the Newton polyhedron}, Math. Annalen {\bf 346} (2010)
 no. 4, 857-890. 
 
\noindent [Gr] P. T. Gressman, {\it Uniform Sublevel Radon-like Inequalities}, J. Geom. Anal. {\bf 23} (2013), no. 2,
611-652.

\noindent [IkKM] I. Ikromov, M. Kempe, and D. M\"uller, {\it Estimates for maximal functions associated
to hypersurfaces in $\R^3$ and related problems of harmonic analysis}, Acta Math. {\bf 204} (2010), no. 2,
151--271.

\noindent [ISa1] A. Iosevich, E. Sawyer, {\it Maximal averages over surfaces},  Adv. Math. {\bf 132} (1997), no. 1, 46-119.

\noindent [ISa2] A. Iosevich and E. Sawyer, {\it Sharp $L^p$ to $L^q$ estimates for a class of averaging operators}, Ann. Inst. Fourier, Grenoble
{\bf 46} (1996), no. 5, 1359-1384.

\noindent [ISaS] A. Iosevich, E. Sawyer, and A. Seeger, {\it On averaging operators associated with convex hypersurfaces of finite type}, 
J. Anal. Math. {\bf 79} (1999), 159-187.

\noindent [Ka1] V. N. Karpushkin, {\it A theorem concerning uniform estimates of oscillatory integrals when
the phase is a function of two variables}, J. Soviet Math. {\bf 35} (1986), 2809-2826.

\noindent [Ka2] V. N. Karpushkin, {\it Uniform estimates of oscillatory integrals with parabolic or 
hyperbolic phases}, J. Soviet Math. {\bf 33} (1986), 1159-1188.

\noindent [L], W. Littman, {\it $L^p$-$L^q$ estimates for singular integral operators arising from hyperbolic equa-
tions}, Partial differential equations. Proc. Sympos. Pure Math. {\bf 23} (1973), 479-481.

\noindent [O] D. M. Oberlin, {\it Convolution with measures on hypersurfaces}, Math. Proc. Camb. Phil. Soc.
{\bf 129} (2000), no. 3, 517-526.

\noindent [R1] B. Randol, {\it On the Fourier transform of the indicator function of a planar set}, Trans. Amer. Math. Soc. {\bf 139} 1969 271-278.

\noindent [R2] B. Randol, {\it On the asymptotic behavior of the Fourier transform of the indicator function of a convex set}, Trans. Amer. Math. Soc. 
{\bf 139} 1969 279-285. 

\noindent [S] A. Seeger, {\it Radon transforms and finite type conditions}, J. Amer. Math. Soc. {\bf 11} (1998), no. 4, 869-897.

\noindent [St] B. Street, {\it Sobolev spaces associated to singular and fractional Radon transforms}, Rev. Mat. Iberoam. {\bf 33} (2017), no. 2, 633-748.

\noindent [Ste] E. M. Stein, {\it $L^p$ boundedness of certain convolution operators}, Bull. Amer. Math. Soc. {\bf 77} (1971), 404-405.

\noindent [Str] R. Strichartz, {\it Convolutions with kernels having singularities on the sphere}, Trans. Amer. Math. Soc. {\bf 148} (1970), 461-471.

\noindent [V] A. N. Varchenko, {\it Newton polyhedra and estimates of oscillatory integrals}, Functional 
Anal. Appl. {\bf 18} (1976), no. 3, 175-196.

\end{document}